\documentclass[11pt,reqno]{amsproc}

\title[Electroconvection models]{On some electroconvection models}

\author[P. Constantin]{Peter Constantin}
\address{Department of Mathematics, Princeton University, Princeton, NJ 08544}
\email{const@math.princeton.edu}

\author[T. Elgindi]{Tarek M. Elgindi}
\address{Department of Mathematics, Princeton University, Princeton, NJ 08544}
\email{tme2@.princeton.edu}

\author[M. Ignatova]{Mihaela Ignatova}
\address{Department of Mathematics, Princeton University, Princeton, NJ 08544}
\email{ignatova@math.princeton.edu}

\author[V. Vicol]{Vlad Vicol}
\address{Department of Mathematics, Princeton University, Princeton, NJ 08544}
\email{vvicol@math.princeton.edu}

\usepackage[margin=1in]{geometry}
\usepackage{amsmath, amsthm, amssymb}
\usepackage{times}
\usepackage{color}
\usepackage{hyperref}
\usepackage{graphicx}
\newcommand{\pa}{\partial}
\newcommand{\la}{\label}
\newcommand{\fr}{\frac}
\newcommand{\na}{\nabla}
\newcommand{\be}{\begin{equation}}
\newcommand{\ee}{\end{equation}}
\newcommand{\ba}{\begin{array}{l}}
\newcommand{\ea}{\end{array}}
\newcommand{\Rr}{{\mathbb R}}

\newcommand{\beg}{\begin}

\newcommand{\D}{\Delta}
\renewcommand{\l}{\Lambda}
\renewcommand{\P}{\mathbb{P}}

\numberwithin{equation}{section}

\begin{document}
\begin{abstract}
{We consider a model of electroconvection motivated by studies of the motion of a two dimensional annular suspended smectic film under the influence of an electric potential maintained at the boundary by two cylindrical electrodes. We prove that this electroconvection model has global in time unique smooth solutions}.  \today.
\end{abstract}
\maketitle
\section{Introduction}
Electroconvection is the flow of fluids and particles driven by electrical forces. There are several studies of electroconvection in the physical literature, pertaining to different types of occurrences of the phenomenon. The interaction of electromagnetic fields with condensed matter is a vast and important subject, with applications ranging from solar magnetohydrodynamics to microfluidics. 
Here we idscuss a particular system, in which a charge distribution interacts with a fluid in a geometrically constrained situation. The fluid is confined to a very thin region, and a voltage difference is maintained by electrodes situated at the boundaries of the region. Physical experiments \cite{exp} and numerical studies \cite{num} consider the flow of an annular suspended smectic film. Despite the non-Newtonian nature of the constituent, the model describes the fluid by Navier-Stokes equations confined to a fixed two dimensional region (an annulus in the cited studies). The Navier-Stokes equations are driven by body forces due to the electrical charge density and the potential. The charge density is transported  by the electric potential and by the flow.
The electric potential is determined by three dimensional equations in the whole space: the physics is inherently nonlocal. 
\vspace{-1cm}
\begin{figure}[htb!]
\begin{center}
\setlength{\unitlength}{1.5in} 
\begin{picture}(2.02083,2)
\put(0.25,0){\includegraphics[height=3in]{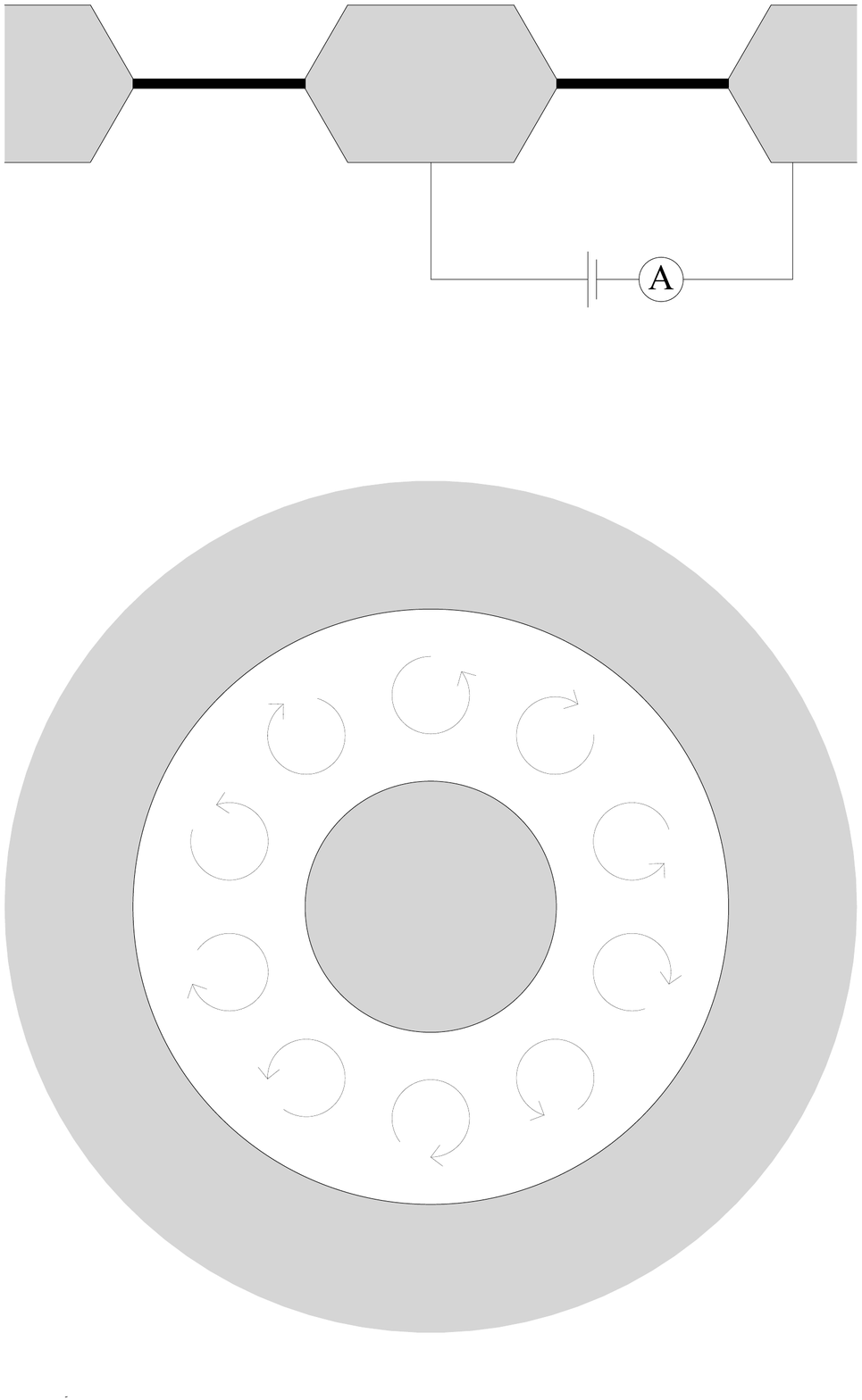}}
\put(0.46,0.6){\mbox{{\tiny $K_1$}}}
\put(1.4,0.6){\mbox{{\tiny $K_1$}}}
\put(0.69,0.6){\mbox{{\tiny $\Omega$}}}
\put(0.94,0.6){\mbox{{\tiny $K_2$}}}
\put(0.94,1.65){\mbox{{\tiny $K_2$}}}
\put(1.22,0.6){\mbox{{\tiny $\Omega$}}}
\put(0.46,1.65){\mbox{{\tiny $K_1$}}}
\put(0.69,1.7){\mbox{{\tiny $\Omega$}}}
\put(1.22,1.7){\mbox{{\tiny $\Omega$}}}
\put(1.4,1.65){\mbox{{\tiny $K_1$}}}
\end{picture}
\caption{Schematic of the experiment considered in~\cite{exp,num}. Side view and top view.}
\end{center}
\end{figure}

The potential obeys 
\[
-\D_3 \Phi = 2q\delta_{\Omega}
\]
with $\Omega\subset \{(x,0)\left|\right. \; x\in \Rr^2\}\cap (\Rr^3\setminus (K_1\cup K_2))$, and $\delta_\Omega$ is the Dirac mass on $\Omega$. The factor $2$ is due to the fact that the film has two sides. The term $2q\delta_{\Omega}$ is the charge density in the limit of zero thickness of the film. Here $\D_3$ is the 3D Laplacian. All the rest of the derivatives below will be 2D. We tacitly identify $\Rr^2$ with $\{(x,0)\left |\right. x\in\Rr^2\}$. The regions $K_1\subset \Rr^3$ and $K_2\subset\Rr^3$ determine the electrodes (disjoint, with smooth or Lipschitz boundaries), and the smooth bounded open set $\Omega\subset \Rr^2\times\{0\}$ shares boundaries with them. The electrodes are maintained at different voltages,
\[
\Phi_{\pa K_1} = V, \quad \Phi_{\pa K_2} = 0,
\]
and we have $\Phi\to 0$ at infinity.

The charge density $q$ evolves in time 
\be
\pa_t q + u\cdot\na q = \D\Phi_{\left |\right. \Omega},
\la{chargeq}
\ee
together with the fluid
\be
\pa_t u + u\cdot\na u - \D u + \na p = -q\na\Phi_{\left |\right. \Omega}.
\la{fluideq}
\ee
The fluid is two dimensional, incompressible
\be
\na\cdot u = 0,
\la{inc}
\ee
and adheres to the boundary, $u_{\left | \right.\pa\Omega} = 0.$ The charge density determines the three dimensional potential
\[
\Phi = \Phi_1 + \Phi_2
\]
where $\Phi_2$ is the 3D Newtonian potential ({fundamental solution of the 3D Laplacian}) of the distribution $2q\delta_{\Omega}$
\[
\Phi_2(x,z) = \fr{1}{4\pi}\int_{\Omega}\fr{2q(y)}{\sqrt{z^2 + (x-y)^2}}dy,
\]
for $x\in\Rr^2$, $z\in\Rr$, and $\Phi_1$ solves
\[
\Delta_3 \Phi_1 = 0 \quad{\mbox{in}}\; \Rr^3\setminus (K_1\cup K_2), 
\]
with boundary conditions
\[
{\Phi_1}_{\left|\right. \;\pa K_1} = V-{\Phi_2}_{\left|\right. \;\pa K_1},\quad
{\Phi_1}_{\left|\right. \;\pa K_2} = -{\Phi_2}_{\left|\right. \;\pa K_2},\quad 
\Phi_1\to 0\; {\mbox{at infinity}}.
\]
The Newtonian potential of the charge distribution, $\Phi_2$, has a jump singularity in the normal derivative
at the points of continuity of $q$ in $x\in\Omega$ 
\[
-[\pa_z \Phi_2]_{\left |\right. z =0} = 2q(x).
\]
{Here we denote the jump of a function $f$ across $z=0$ by $[f]_{\left |\right. z =0}  = \lim_{z\downarrow 0} f(z) - \lim_{z \uparrow 0} f(z)$.}
Thus, the Newtonian potential of the charge $2q\delta_{\Omega}$ with $q$ continuous in $\Omega$ can be seen as the distributional solution of the two phase problem
\[
\D_3 \Phi_2 = 0, \quad {\mbox{in}} \;  z\neq 0,
\]
\[
[\Phi_2]_{\left |\right. z =0} = 0,\quad x\notin{\pa \Omega}
\]
\[
-[\pa_z\Phi_2]_{\left |\right. z =0} = \left\{\ba 
2q(x),\quad {\mbox{for}} \; x\in \Omega\\
0\; \quad \quad {\mbox{for}} \; x\notin {\overline{\Omega}}.
\ea
\right.
\]
Note that $-\pa_z\Phi_2$ is the Poisson integral of $q\chi_{\Omega}$ where $\chi_{\Omega}$ is the characteristic function of the set $\Omega$. The Poisson integral equals
\[
-\pa_z\Phi_2(x,z) = e^{-z\l_{\Rr^2}}(q\chi_{\Omega})(x)
\]
for $z>0$, and therefore
\[
\D\Phi_2(x,z) = \pa_ze^{-z\l_{\Rr^2}}(q\chi_{\Omega})(x) = -\l_{\Rr^2}\left(e^{-z\l_{\Rr^2}}(q\chi_{\Omega})\right)(x)
 \]
holds for $z>0$. Passing to the limit $z\downarrow 0$ we obtain that the contribution due to $\Phi_2$  to the right side of \eqref{chargeq} is given by
\[
\D\Phi_2(x,0) = -\l_{\Rr^2}^*(q\chi_{\Omega})
\]
where $\l_{\Rr^2}: H^1(\Rr^2)\to L^2(\Rr^2)$ is the square root of the Laplacian $\l_{\Rr^2} = \sqrt{-\D}$, and $\l_{\Rr^2}^*: L^2(\Rr^2)\to H^{-1}(\Rr^2)$ is its adjoint. The equality $\l^*_{\Rr^2}(q\chi_{\Omega}) = \l_{\Rr^2}(q\chi_{\Omega})$ is valid only if $q\chi_{\Omega}\in H^1(\Rr^2)$. On the other hand,  a direct calculation shows if $q_{\left |\right. \pa \Omega} = 0$ that
\[
\D\Phi_2(x,z) = -\fr{1}{2\pi}\int_{\Rr^2}\fr{x_i-y_i}{(z^2 + |x-y|^2)^{\fr{3}{2}}}\chi_{\Omega}(y)\fr{\pa q(y)}{\pa y_i} dy,
\]
and passing to the limit $z\downarrow 0$ we obtain
\[
\D\Phi_2(x,0) = (R_i(\chi_{\Omega}\pa_i q))(x)
\]
if $\na q\in L^2(\Omega)$. Here $R_i$ are Riesz transforms in $L^2(\Rr^2)$. Note that the condition $q\in H_0^1(\Omega)$ is exactly the condition required for $\chi_{\Omega} q\in H^1(\Rr^2)$, if $\pa\Omega$ is smooth. It is easy to see that
\[
\pa_iR_i(\chi_{\Omega}q) = -\l_{\Rr^2}^*(\chi_{\Omega}q)
\]
holds in general, if $q\in L^2(\Omega)$. 

Regarding $\Phi_1$ let us note that because $\Phi_2(x,z)$ belongs to $C^{\alpha}(\Rr^3)$ with $0<\alpha<1$ (it is actually Lipschitz) if $q$ is bounded, it follows that $\Phi_1\in C^{2,\alpha}(\Rr^3\setminus (K_1\cup K_2)$.

Unfortunately, the condition that 
$q\in H_0^1(\Omega)$ is not maintained under the evolution, even if $u=0$ and $\Phi_1 =0$.
In that case the solution of
\[
\pa_tq = -\l_{\Rr^2}^*(\chi_{\Omega}q)
\]
is given by the Poisson kernel,
\[
q(x,t) = \fr{1}{2\pi}\int_{\Omega}\fr{t}{(t^2 + |x-y|^2)^{\fr{3}{2}}}q_0(y)dy.
\]
which manifestly does not keep $q_{\left |\right. \pa\Omega}= 0.$

The addition of the contribution from the potential $\Phi_1$ does not help matters. Singularities in $q$ appear at the shared boundaries, and their analysis is not our priority here. The main message of the calculation above is that nonlocal dissipation appears naturally when electrical (or electro-magnetic) fields interact with confined condensed matter. 

\vspace{-0.5cm}
\begin{figure}[htb!]
\begin{center}
\setlength{\unitlength}{1.5in} 
\begin{picture}(2.02083,2)
\put(0.25,0){\includegraphics[height=3in]{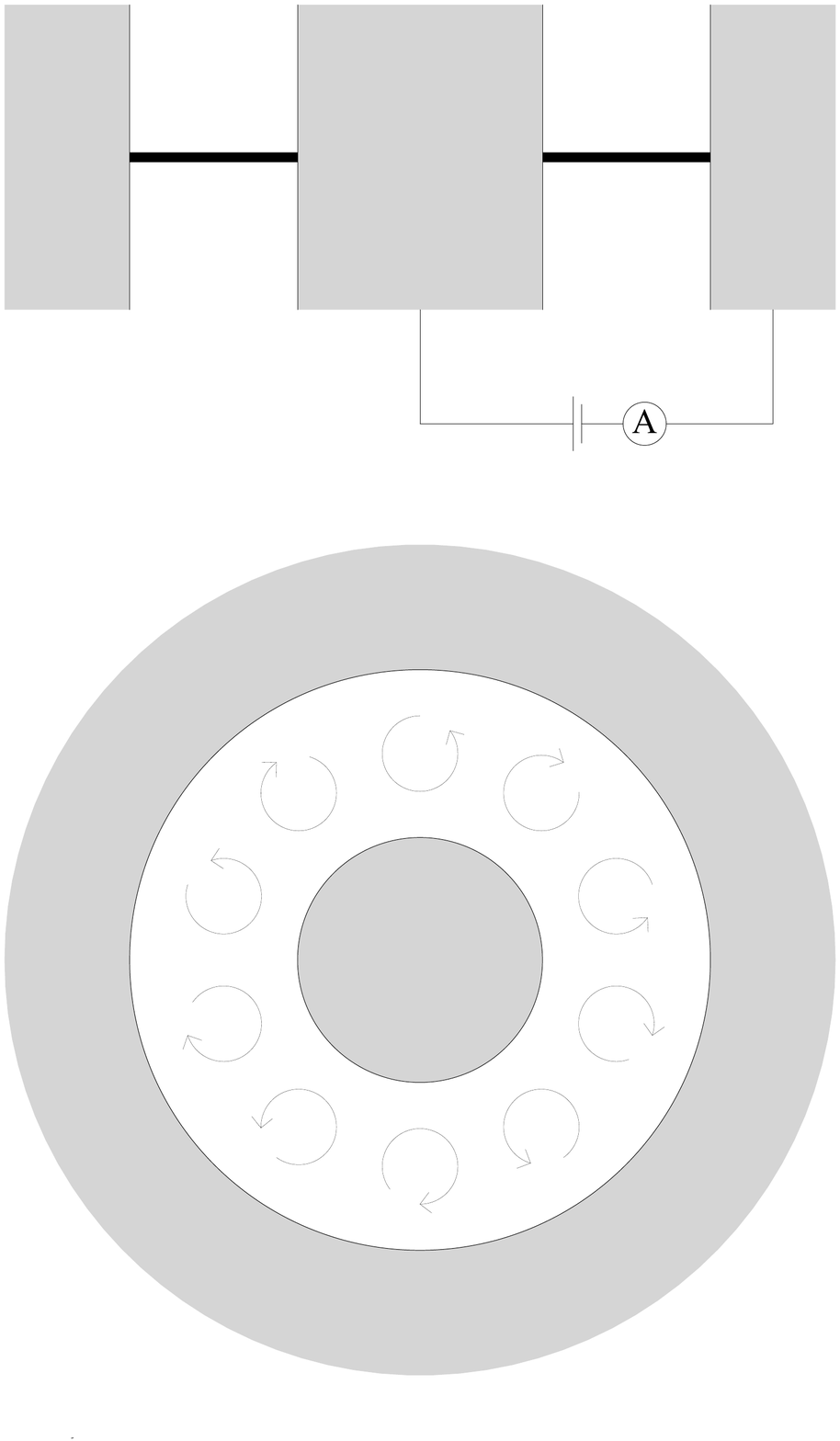}}
\put(0.46,0.6){\mbox{{\tiny $K_1$}}}
\put(1.4,0.6){\mbox{{\tiny $K_1$}}}
\put(0.69,0.6){\mbox{{\tiny $\Omega$}}}
\put(0.94,0.6){\mbox{{\tiny $K_2$}}}
\put(0.94,1.65){\mbox{{\tiny $K_2$}}}
\put(1.22,0.6){\mbox{{\tiny $\Omega$}}}
\put(0.46,1.65){\mbox{{\tiny $K_1$}}}
\put(0.69,1.7){\mbox{{\tiny $\Omega$}}}
\put(1.22,1.7){\mbox{{\tiny $\Omega$}}}
\put(1.4,1.65){\mbox{{\tiny $K_1$}}}
\end{picture}
\caption{Schematic of the model considered in this paper. Side view and top view.}
\end{center}
\end{figure}

In this paper we avoid the boundary singularities by considering the electrods $K_1$ and $K_2$ to extend to infinity in the $z$ direction. This way the potential $\Phi_1$ {is} $z$-independent, and does not contribute to {the $q$ evolution \eqref{chargeq}}. To be more specific, we consider a connected open domain $\Omega$ with smooth boundary in $\Rr^2$. The domain is not simply connected, and $\pa\Omega = \Gamma_1 \cup \Gamma_2$ with
$\Gamma_1\cap \Gamma_2 = \emptyset$ and ${\rm dist}(\Gamma_1,\Gamma_2)>0$. (An annular region is such a domain, but that double connectivity is not important). We solve the problem
\[
-\D_3\Phi_2(x,z) = 2q(x)\delta_{\Omega}
\]
for $x\in\Omega$ together with the boundary condition ${\Phi_2}_{\left|\right. \pa\Omega} = 0$ by setting
\[
\Phi_2(x,z) =  
\begin{cases}
(e^{-z\l_D}\l_D^{-1}q)(x), \quad {\mbox{if}} \; z>0,\\  
(e^{z\l_D}\l_{D}^{-1}q)(x), \; \; \quad {\mbox{if}} \; z<0,
\end{cases}
\]
where $\l_D$ is the square root of the Dirichlet Laplacian.  
If we take $K_1 = \Gamma_1\times \Rr$ and $K_2 = \Gamma_2\times\Rr$, i.e., we consider vertical electrods at the boundaries of the domain, then the harmonic function $\Phi_1(x)$, solving $\D \Phi_1 = 0$ in $\Omega$ and ${\Phi_1}_{\left|\right. \Gamma_1}= V$, ${\Phi_1}_{\left |\right. \Gamma_2} = 0$ solves also
$\D_3\Phi_1 = 0$ with boundary conditions ${\Phi_1}_{\left|\right. K_1}= V$, ${\Phi_1}_{\left |\right. K_2} = 0$ in $\Omega\times \Rr$. The function $\Phi_1$ is $z$-independent and harmonic in $x$. Because $\Phi_1$ is smooth, it might influence the stability of solutions to \eqref{chargeq}--\eqref{fluideq}, but not their regularity. For simplicity of exposition we take $\Phi_1 =0$. 

The system (\ref{chargeq}), (\ref{fluideq}), (\ref{inc})  becomes:
\be
\left\{
\ba
\pa_t q + u\cdot\na q + \l_D q = 0,\\
\pa_t u + u\cdot\na u -\D u +\na p = {-q \na\l_D^{-1}q},\qquad  \mbox{in }\Omega \\
\na\cdot u = 0,
\ea
\right.
\la{sys}
\ee
with homogeneous Dirichlet boundary conditions for $u$ on $\partial \Omega$. In the rest of the paper we study the  regularity of solutions to \eqref{sys}. Our main result, Theorem~\ref{main} below, shows that if $u_0 \in H_0^1(\Omega) \cap H^2(\Omega) \cap \P(L^2(\Omega))$ and $q_0 \in H_0^1(\Omega) \cap H^2(\Omega)$, then solutions to \eqref{sys} exist for all time, are smooth and uniquely determined by initial data.

\section{Preliminaries}
We consider the Dirichlet Laplacian in a bounded open domain
$\Omega\subset \Rr^d$ with smooth boundary. We denote by $\D$ the Laplacian operator with homogeneous Dirichlet boundary conditions. Its $L^2(\Omega)$ - normalized eigenfunctions are denoted $\phi_j$, and its eigenvalues counted with their multiplicities are denoted $\mu_j$: 
\be
-\D \phi_j = \mu_j \phi_j.
\la{ef}
\ee
It is well known that $0<\mu_1\le...\le \mu_j\to \infty$  and that $-\D$ is a positive selfadjoint operator in $L^2(\Omega)$ with domain ${\mathcal{D}}\left(-\D\right) = H^2(\Omega)\cap H_0^1(\Omega)$.
Functional calculus is defined using the eigenfunction expansion. In particular
\be
\left(-\D\right)^{\alpha}f = \sum_{j=1}^{\infty}\mu_j^{\alpha} f_j \phi_j
\la{funct}
\ee
with 
\[
f_j =\int_{\Omega}f(y)\phi_j(y)dy
\]
for $f\in{\mathcal{D}}\left(\left (-\D\right)^{\alpha}\right) = \{f\left |\right. \; (\mu_j^{\alpha}f_j)\in \ell^2(\mathbb N)\}$.
We denote by
\be
\l^s = \left(-\D\right)^{\fr{s}{2}}
\la{lambdas}
\ee
the fractional powers of the Dirichlet Laplacian and by 
\be
\|f\|_{s,D}^2 = \sum_{j=1}^{\infty}\mu_j^{s}f_j^2.
\la{norms}
\ee
the norms in  ${\mathcal{D}}\left (\l^s\right)$. It is well-known that
\[
{\mathcal{D}}\left( \l \right) = H_0^1(\Omega).
\]

We recall the C\'{o}rdoba-C\'{o}rdoba inequality \cite{cc} for bounded domains, {established in} \cite{ci}:
\beg{prop}{\la{cordoba}}
Let $\Phi$ be a $C^2$ convex function satisfying $\Phi(0)= 0$. Let $f\in C_0^{\infty}(\Omega)$ and let $0\le s\le 2$. Then
\be
\Phi'(f)\l^s f - \l^s(\Phi(f))\ge 0
\la{cor}
 \ee
holds pointwise almost everywhere in $\Omega$.
\end{prop}
We use also the following commutator estimate {proven in} \cite{ci}:
\beg{thm}{\la{comthm}}
Let a vector field $v$ have components in $B(\Omega)$ where $B(\Omega) = W^{2,d}(\Omega)\cap W^{1,\infty}(\Omega)$, if $d\ge 3$, and $B(\Omega) = W^{2,p}(\Omega)$ with $p>2$, if $d=2$. Assume that the normal component of
 the trace of $v$ on the boundary vanishes, 
\[
v_{\left |\right. \pa \Omega}\cdot n = 0
\]
(i.e the vector field is tangent to the boundary). 
There exists a constant $C$ such that
\be
\|[v\cdot\na ,\l]f\|_{\fr{1}{2}, D} \le C\|v\|_{B(\Omega)}\|f\|_{\fr{3}{2},D}
\la{comtwo}
\ee
holds for any $f$ such that $f\in {\mathcal{D}}\left(\l^{\fr{3}{2}}\right)$,
where
\[
\|v\|_{B(\Omega)} = \|v\|_{W^{2,d}(\Omega)} + \|v\|_{W^{1,\infty}(\Omega)}
\]
if $d\ge 3$ and
\[
\|v\|_{B(\Omega)} = \|v\|_{W^{2,p}(\Omega)}
\]
with $p>2$, if $d=2$.
\end{thm}
This result is proved in \cite{ci} using the method of harmonic extension.  It is used to prove an existence theorem for linear equations of transport and nonlocal diffusion \cite{ci}:

\beg{thm}{\la{linthm}} Let $u\in L^2(0,T; B(\Omega)^d)$ be a vector field parallel to the boundary. Then the equation 
\be
\pa_t\theta + u\cdot\na\theta + \l\theta = 0
\la{thetaeq}
\ee
with initial data $\theta_0\in H_0^1(\Omega)\cap H^2(\Omega)$ has unique solutions belonging to
\[
\theta\in L^{\infty}(0,T; H^2(\Omega)\cap H_0^1(\Omega))\cap L^2(0,T; H^{2.5}(\Omega)).
\]
If the initial data $\theta_0\in L^p(\Omega)$, $1\le p\le \infty$, then
\be
\sup_{0\le t \le T}\|\theta(\cdot, t)\|_{L^p(\Omega)} \le \|\theta_0\|_{L^p(\Omega)}
\la{lp}
\ee
holds.
\end{thm}
We need also the fact that, for $d=2$,
\be
\|f\|_{L^4(\Omega)}\le C\|f\|_{\fr{1}{2}, D}.
\la{lfour}
\ee 
Background material and applications of the method of harmonic extension can be found  in \cite{cabre}. 

We recall now some basic notions concerning the Navier-Stokes equations \cite{cfbook}.
The Stokes operator 
\be
A = -\P\D\la{stokes}
\ee
is defined via the Leray-Hodge projector 
\be
\P: L^2(\Omega)^d\to H=\{u\left |\right. \; u\in L^2(\Omega)^d, \; \na\cdot u = 0\}
\la{lerayp}
\ee
The domain of $A$ in $L^2$ is ${\mathcal D}(A) = H_0^1(\Omega)^d\cap H^2(\Omega)^d\cap H$. The operator is positive, 
\[
(Au,u)_H = \int_{\Omega}|\na u|^2dx, \quad \forall \; u\in {\mathcal{D}}(A),
\]
elliptic and injective,
\be
\|u\|_{H^2(\Omega)^d}\le C \|Au\|_{H},
\la{ell}
\ee
and its inverse $A^{-1}$ is compact. Functional calculus is defined using eigenfunction expansion. The eigenvalues of $A$ are denoted $\lambda_j$, 
$0<\lambda_1\le \dots \lambda_j\le \dots \to\infty$, the eigenfunctions $w_j\in{\mathcal {D}}(A)$,
\[
Aw_j = \lambda_j w_j.
\]
The square root $A^{\fr{1}{2}}$ satisfies
\be
\|A^{\fr{1}{2}}v\|_{H} = \|\na v\|_{L^2(\Omega)}
\la{katoa}
\ee
for any $v\in H\cap H_0^1(\Omega)^d$. The nonlinear term in the Navier-Stokes equations is
\be
B(u,u) = \P(u\cdot\na u).
\la{buu}
\ee
It has the property that
\[
(B(u,u), u)_H = 0, 
\]
for all $u\in H\cap H_0^1(\Omega)^d$, $d\le 3$.
In addition, for $d=2$, using
\be
\|u\|_{L^4(\Omega)}\le C\|u\|_{L^2(\Omega)}^{\fr{1}{2}}\|\na u\|_{L^2(\Omega)}^{\fr{1}{2}}
\la{ulfour}
\ee
and
\be
\|\na u\|_{L^4(\Omega)}\le C\|\na u\|_{L^2(\Omega)}^{\fr{1}{2}}\| u\|_{H^2(\Omega)}^{\fr{1}{2}}
\la{naulfour}
\ee
and the ellipticity of the Stokes operator, we obtain
\be
\|B(u,u)\|_{H}\le C \|u\|_{L^2(\Omega)}^{\fr{1}{2}}\|\na u\|_{L^2(\Omega)}\|Au\|_{H}^{\fr{1}{2}},
\la{buuau}
\ee
valid for all $u\in{\mathcal{D}}(A)$.
Also, using elliptic regularity we have
\be
\|A^{\fr{1}{2}}B(u,u)\|_H \le C\left[ \|u\|_H^{\fr{1}{2}}\|Au\|^{\fr{3}{2}}_H + \|\na u\|_{L^2(\Omega)}\|Au\|_H\right]
\la{abuu}
\ee
for $u\in{\mathcal{D}}(A)$.
Indeed, the right hand side bounds $\na(u\cdot\na u)= u\cdot\na\na u + \na u \na u $ in $L^2(\Omega)$, and, because $B(u,u) = u\cdot\na u + \na \pi$ with
$-\D \pi = \na(u\cdot\na u)$, with homogeneous Neumann boundary condition $\fr{\pa \pi}{\pa n}= 0$ at $\pa\Omega$, it follows by elliptic regularity that $\na \na \pi$ obeys the same $L^2(\Omega)$ bound.  
\section{A base model}

We consider the system formed by the incompressible Navier-Stokes equations 
\be
\pa_t u + u\cdot\na u -\Delta u + \na p = -qRq 
\la{nse}
\ee
with 
\be
\na\cdot u = 0,
\la{divu}
\ee
and with
\be
Rq = \na\l^{-1}q,
\la{rieszd}
\ee
coupled with the evolution of the charge density
\be
\pa_t q + u\cdot\na q +\l q = 0.
\la{qeq}
\ee
The equations hold for $x\in\Omega$, $t\ge 0$, and the initial data $u_0$ and $q_0$ are smooth.
Our main {result} is:
\beg{thm}\la{main} Let $u_0\in{\mathcal{D}}(A)$,  $q_0\in H_0^1(\Omega)\cap H^2(\Omega)$ and let $T$ be arbitrary. Then, there exist and are unique solutions of the problem (\ref{nse}), (\ref{divu}), (\ref{qeq}) in $d=2$ on the time interval $[0,T]$ obeying
\be
u\in L^{\infty}(0,T; {\mathcal {D}}(A)) \cap L^2(0,T; {\mathcal {D}}(A^{\fr{3}{2}}))
\la{ub}
\ee
and
\be
q\in L^{\infty}(0,T; W^{1,4}(\Omega)) \cap L^2(0,T; H^2(\Omega))
\la{qb}
\ee
with explicit bounds that depend only on the initial data and not on $T$. More precisely, there exists an explicit function of one variable, $C[N]$, with double exponential growth in $N$, such that
\begin{align}
&\| u\|_{ L^{\infty}(0,T; {\mathcal {D}}(A)) \cap L^2(0,T; {\mathcal {D}}(A^{\fr{3}{2}}))} + \|q\|_{ L^{\infty}(0,T; W^{1,4}(\Omega)) \cap L^2(0,T; H^2(\Omega))} \notag \\
& \le C\left [\|u_0\|_{{\mathcal{D}}(A)} + \|q_0\|_{{\mathcal {D}}(-\D)}\right].
\la{fin}
\end{align}
\end{thm} 

In order to prove the theorem, we construct solutions by an approximation procedure {(see \eqref{pmnse}, \eqref{pmqeq}, \eqref{idum} below)}, prove a priori estimates on the approximants {( see \eqref{qhone}, \eqref{super}, \eqref{qwonefour} below)}, and pass to the limit {via the Aubin-Lions compactness theorem}.

We consider Galerkin approximations for $u$. These are defined using the projectors $\P_m$:
\be
\P_m u = \sum_{j=1}^m (u,w_j)_H\,w_j
\la{pmu}
\ee
The approximate system is
\be
\pa_t u_m + Au_m +\P_m(B(u_m,u_m)) = \P_m(qRq)
\la{pmnse}
\ee
for $u_m\in \P_m H$,
coupled with
\be
\pa_t q + u_m\cdot\na q + \l q =0
\la{pmqeq}
\ee
with initial data $q_0\in H_0^1(\Omega)\cap H^2(\Omega)$. The initial data for $u_m$ are
\be
u_m(0) = \P_mu_0
\la{idum}
\ee
where $u_0\in{\mathcal{D}}(A)$ is the initial velocity in our problem. The system is thus a system of nonlinear ordinary differential equations (\ref{pmnse}) coupled to a linear
transport and nonlocal diffusion partial differential equation (\ref{pmqeq}). {By Theorem \ref{linthm} (see \cite{ci}) we have that} the linear equation (\ref{pmqeq}) has unique solutions in
\[
q\in L^{\infty}(0,T; H^2(\Omega)\cap H_0^1(\Omega))\cap L^2(0,T; H^{2.5}(\Omega))
\]
if $q_0\in  H_0^1(\Omega)\cap H^2(\Omega)$, as long as $u_m\in L^2(0,T; B(\Omega)^d)$ is a vector field wich is parallel to the boundary. Because $u_m$ is a finite linear combination of eigenfunctions of the Stokes operator, smooth and vanishing at the boundary, the only issue is whether some norm of $u_m$ stays finite and square integrable in time. The reason we chose the full PDE for $q$ rather than some approximation is so that we can use $L^{\infty}$ bounds.
We prove this in conjunction with a priori estimates on $q$ which follow from Proposition \ref{cordoba}:
\be
\sup_{0\le t\le T}\|q\|_{L^p(\Omega)}\le \|q_0\|_{L^p(\Omega)}.\la{lpq}
\ee
These are valid for any $p$, $1\le p \le \infty$. These inequalities are justified in our situation by the energy bounds for $u_m$ which follow below. Because $R:L^2(\Omega)\to L^2(\Omega)$ are bounded, it follows then that the forcing $qRq$ in the right hand side of (\ref{pmnse}) is bounded uniformly in $m$
\be
\|\P_m(qRq)\|_{L^{\infty}(0,T; L^2(\Omega))}\le C\|q_0\|_{L^{\infty}(\Omega)}\|q_0\|_{L^2(\Omega)}.
\la{l2qrq}
\ee
The energy inequality
\be
\fr{d}{dt}\|u_m\|^2_{H} + \|\na u_m\|^2_{L^2(\Omega)} \le C\|qRq\|_{L^2(\Omega)}^2 
\la{enum}
\ee
implies that $u_m\in L^{\infty}(0,T; L^2(\Omega)^2)\cap L^2(0,T; H^1(\Omega)^2)$ are uniformly bounded in $m$:
\be
\sup_{0\le t\le T}\|u_m\|_H^2 + \int_0^T\|\na u_m\|^2_{L^2(\Omega)} \le \|u_0\|^2_H + C \|q_0\|_{L^{\infty}(\Omega)}^2\|q_0\|_{L^2(\Omega)}^2.
\la{enb}
\ee
Strictly speaking, because the system formed by (\ref{pmnse}) and (\ref{pmqeq}) is nonlinear, we have to redo the proof of Theorem \ref{linthm} for this nonlinear case. This involves taking an additional approximation, at fixed $m$, of (\ref{pmqeq}) by eigenfunction expansions of the Laplacian, proving energy bounds for it, and passing to the limit. This approximation does not preserve {the $L^{\infty}$ norm}, and because of that, passage to the limit is done first for short time. The limit equation is the nonlinear system (\ref{pmnse}), (\ref{pmqeq}), and because of smoothness of $u_m$ we obtain the uniform bounds (\ref{lpq}) and (\ref{enb}) a posteriori. Then we extend the solution. A fixed point argument using a semigroup method needs to be avoided because the system is not semilinear. 
Here are a few details of this procedure: We couple 
\be
\pa_t u_m + Au_m +\P_m(B(u_m,u_m)) = \P_m(q_nRq_n)
\la{pmnnse}
\ee
for $u_m\in \P_m H$, with initial data $u_m(0) = \P_mu_0$, to
\be
\pa_t q_n + P_n(u_m\cdot\na q_n) + \l q_n =0
\la{pnmq}
\ee
with 
\be
P_n (f) = \sum_{j=1}^n (f,\phi_j)_{L^2(\Omega)}\,\phi_j
\la{pn}
\ee
and initial data $q_n(0) = P_nq_0$. Thus (\ref{pmnnse})--(\ref{pnmq}) is a system of ODEs which have solutions on a maximal time interval. How long  this time interval {is} depends on bounds. We use energy bounds employing the commutator estimate (\ref{comtwo}) and the proof of Theorem \ref{linthm} to obtain a priori estimates. The basic estimate concerns the $H^2$ norms. 
We apply $\l$ to (\ref{pnmq}), and use the commutator:
\be
\pa_t\l q_n + \l(\l q_n) + P_n(u_m\cdot\na \l q_n) + P_n[\l, u_m\cdot\na] q_n = 0.
\la{lqneq}
\ee
We take the scalar product with $\l^3q_n$:
\begin{align*}
&\int_{\Omega}(u_m\cdot\na\l q_n)\l^3 q_n dx 
= \int_{\Omega}\l^2(u_m\cdot\na\l q_n)\l q_n dx\\
&=\int_{\Omega}\left [(-\Delta u_m)\cdot\na \l q_n - 2\na u_m\cdot\na\na\l q_n\right]\l q_n dx  + \int_{\Omega}(u_m\cdot\na\l^3q_n)\l q_n dx\\
&= \int_{\Omega}\left [(-\Delta u_m)\cdot\na \l q_n - 2\na u_m\cdot\na\na\l q_n\right]\l q_ndx - \int_{\Omega}\l^3q_n(u_m\cdot\na \l q_n)dx\\
&=\int_{\Omega}\left [((-\Delta u_m)\cdot\na \l q_n)\l q_n + 2\na u_m\na\l q_n\na\l q_n\right]dx-\int_{\Omega}(u_m\cdot\na\l q_n)\l^3 q_ndx.
\end{align*}
In the first integration by parts we used the fact that $\l q_n$ is a finite linear combination of eigenfunctions which vanish at the boundary. Then we use the fact that $\l^2 = -\Delta$ is local. In the last equality we integrated by parts using the fact that $\l q_n$ is a finite linear combination of eigenfunctions which vanish at the boundary and the fact that $u_m$ is divergence-free.
It follows that
\[
\int_{\Omega}(u_m\cdot\na\l q_n)\l^3 q_n dx = \fr{1}{2}\int_{\Omega}\left [((-\Delta u_m)\cdot\na \l q_n)\l q_n + 2\na u_m\na\l q_n \na\l q_n\right]dx
\]
and consequently
\[
\left| \int_{\Omega}(u_m\cdot\na\l q_n)\l^3 q_n dx\right| \le C\|u_m\|_{B(\Omega)}\|\l^2 q_n\|_{L^2(\Omega)}^2.
\]
We obtain thus
\begin{align}
\sup_{0\le t\le \tau}\|\l^2 q_n(\cdot, t)\|^2_{L^2(\Omega)} + \int_0^{\tau}\|\l^{\fr{5}{2}}q_n\|^2_{L^2(\Omega)}dt 
&\le C\|\l^2q_0\|^2_{L^2(\Omega)}e^{C\int_0^{\tau}(\|u_m\|_{B(\Omega)}  + \|u_m\|_{B(\Omega)}^2)dt} \notag \\
&\le C\|\l^2q_0\|^2_{L^2(\Omega)}e^{C_m\int_0^{\tau}(\|u_m\|_{L^2(\Omega)} + \|u_m\|_{L^2(\Omega)}^2)dt}
\la{qnb}
\end{align}
on the time of existence interval $[0,\tau)$. We used here the finiteness of $m$ and the smoothness of the eigenfunctions of the Stokes operator:
\be
\|u_m\|_{B(\Omega)}\le C_m\|u_m\|_{L^2(\Omega)}.
\la{umbum}
\ee
The terms with $\|u_m\|^2_{B(\Omega)}$ in the exponent originate from the commutator estimate
(\ref{comtwo}):
\[
\|[\l, u_m\cdot\na ]q_n\|_{\fr{1}{2}, D}\le C\|u_m\|_{B(\Omega)}\|q_n\|_{\fr{3}{2}, D}
\]
which are followed by a Young inequality and use of the dissipative term $\|\l^{2.5}q_n\|^2_{L^2(\Omega)}$.
On the other hand, the basic energy estimate for (\ref{pmnnse}) gives
\be
\fr{d}{dt}\|u_m\|^2_H + \|\na u_m\|^2_{L^2(\Omega)} \le C\|q_nRq_n\|^2_{L^2(\Omega)}
\la{enumn}
\ee
and therefore we deduce
\be
\sup_{0\le t\le \tau}\|u_m\|^2_H \le \|u_0\|^2_H + C\int_0^{\tau}\|q_nRq_n\|^2_{L^2(\Omega)}.
\la{enumnb}
\ee
In $d=2$ we have $H^2(\Omega)\subset L^{\infty}(\Omega)$ and we can replace the $L^2$ norm in the right hand side of the inequality (\ref{enumnb}) above by the fourth power of the $H^2$ norm of $q_n$, 
resulting in 
\[
\sup_{0\le t \le \tau} \|u_m\|^2_H \le C\tau\sup_{0\le t\le \tau}\|q_n\|_{H^2(\Omega)}^4 + \|u_0\|^2_H.
\]
Combining this with (\ref{qnb}) we can close the estimates and show that time of existence for the ODE (\ref{pmnnse}), (\ref{pnmq}) is independent of $n$.
Passing to the limit $n\to\infty$ at fixed $m$, we obtain solutions $(q,u_m)$ of (\ref{pmnse}), (\ref{pmqeq}) on  a short time interval $[0,T_m)$ with $T_m>0$, depending on $m$ and norms of initial data. On this short time interval we obtain a posteriori the uniform $L^{\infty}$ bounds (\ref{lpq}) and the uniform energy bound (\ref{enb}).  With a similar calculation
of the $H^2$ norm of $q$ with this information and obtain that the bound on $\|q\|_{H^2(\Omega)}$ depends only on initial data (at $t=0$) and on $m$. More precisely, the $H^2$ norm of $q$ obeys the same inequality as (\ref{qnb}) and thus it grows at most exponentially in time
\[
\|q(t)\|_{H^2(\Omega)}^2 \le e^{\gamma_m t}\|q_0\|_{H^2(\Omega)}^2
\] 
with $\gamma_m$ depending on $m$ and proportional to the right hand side of (\ref{enb}).
Higher derivatives can be bounded as well. If $T_m <T$, the arbitrary time we are considering, then we can uniquely extend the solution beyond $T_m$ by repeating the procedure, and obtain again the bound above. This proves the global existence of the solutions of the approximate system. 

Using classical methods for $2D$ Navier-Stokes equations it follows next that
$u_m$ are uniformly in $m$ bounded 
\be
u_m\in L^{\infty}(0,T, H^1_{0}(\Omega)^2)\cap L^2(0,T; H^2(\Omega)^2).
\la{ustrong}
\ee
Indeed, we take the scalar product of (\ref{pmnse}) in $H$ with $Au_m$,
and then, using (\ref{buuau}) and Young's inequality we obtain the evolution inequality
\be
\fr{d}{dt}\|\na u_m\|^2_{L^2(\Omega)} + \|Au_m\|^2_{L^2(\Omega)} \le C\|qRq\|_{L^2(\Omega)}^2 + C\|u_m\|_{H}^2\|\na u_m\|_{L^2(\Omega)}^4
\la{ensum}
\ee
and thus, using a Gronwall inequality and the bounds from the energy inequality (\ref{enb}), we obtain (\ref{ustrong}) with bounds independent of $m$. The bounds depend only on the norms of initial data do not depend on $T$.
Note that at this stage we do not yet have uniform bounds in $m$ for $u_m \in L^2(0,T; B(\Omega))$ although obviously we do have each $u_m\in L^2(0,T; B(\Omega))$ in view of the weaker bounds and the fact that $u_m$ are functions in $\P_m H$. 
We take now the scalar product of the equation (\ref{pmqeq}) with $-\D q$. 
\[
\fr{1}{2}\fr{d}{dt}\|\na q\|_{L^2(\Omega)}^2 + \|q\|^2_{\fr{3}{2}, D} =
\int_{\Omega} (u_m\cdot\na q)\D qdx.
\]
We integrate by parts, and  using the facts that $u_m$ vanish on the boundary and are divergence-free, we obtain
\[
\int_{\Omega} (u_m\cdot\na q)\D qdx = -\int_{\Omega}\na q (\na u_m)\na q  dx.
\]
We use a H\"{older} inequality to bound
\[
\left |\int_{\Omega}\na q (\na u_m)\na q  dx\right |\le \|\na q\|_{L^2(\Omega)}\|\na q\|_{L^4(\Omega)}\|\na u\|_{L^4(\Omega)}.
\]
Now, we claim that
\[
\|\na q\|_{L^4(\Omega)}\le C\|q\|_{\fr{3}{2}, D}.
\]
This is true because $R:L^4(\Omega)\to L^4(\Omega)$ is bounded \cite{jerisonkenig}, and because $\l q\in {\mathcal{D}}(\l^{\fr{1}{2}})$; in fact $q\in {\mathcal{D}}(\D)$ by \cite{ci}. (We remark here that since $\partial \Omega$ is smooth, the boundedness of the associated Riesz transforms $R = \nabla (-\Delta_D)^{-1/2}$ follows by a classical argument:
flattening of the boundary, maximal elliptic $L^p$ regularity of the operator $L = {\rm div} A \nabla$ when the matrix $A$ is smooth and uniformly elliptic, and complex interpolation. This argument carries over to the case of Lipschitz domains, albeit only for a restricted range of $p$s, see~\cite{jerisonkenig} for details.) Thus, $\na q = R\l q$ is bounded in $L^4$  using (\ref{lfour}).  Using (\ref{ell}) and (\ref{naulfour}),  we deduce
\[
\left |\int_{\Omega}\na q (\na u_m)\na q  dx\right |\le C\|\na q\|_{L^2(\Omega)}\|q\|_{\fr{3}{2}, D}\|\na u_m\|_{L^2(\Omega)}^{\fr{1}{2}}\|Au_m\|_{H}^{\fr{1}{2}}.
\]
Consequently, after a Young inequality, because of (\ref{ustrong}), it follows that
\be
q\in L^{\infty}(0,T; H^1(\Omega))\cap L^2(0,T; {\mathcal D}(\l^{\fr{3}{2}}))
\la{qhone}
\ee
is bounded a priori, independent on time $T$, in terms only of initial data. 

Now we take (\ref{pmnse}), apply $A$ and take the scalar product with $Au_m$. We obtain the differential inequality
\[
\fr{1}{2}\fr{d}{dt}\|Au_m\|^2_H + \|A^{\fr{3}{2}} u_m\|^2_H \le C\|A^{\fr{1}{2}}\P_m(qRq)\|_H^2 + C\|A^{\fr{1}{2}}\P_m B(u_m,u_m)\|^2_H.
\]
Now
\[
\|A^{\fr{1}{2}}\P_m(qRq)\|_H^2 \le \|A^{\fr{1}{2}}(qRq)\|_H^2 =
\|\na(qRq)\|_{L^2(\Omega)}^2
\]
and
\[
\|\na(qRq)\|_{L^2(\Omega)}^2 \le \|Rq\|_{L^4(\Omega)}^2\|\na q\|_{L^4(\Omega)}^2 + \|q\|_{L^{\infty}(\Omega)}^2\|\na q\|_{L^2(\Omega)}^2 \in L^1((0,T))
\]
because {Riesz transforms are bounded in $L^4(\Omega)$ (cf.~\cite{jerisonkenig})} and because of the bounds (\ref{lpq}) and (\ref{qhone}).  The other term obeys, using (\ref{abuu})
\[
\|A^{\fr{1}{2}}\P_m B(u_m,u_m)\|^2_H \le C \left[\|u_m\|_H\|Au_m\|_H^3 +
\|\na u_m\|_H\|Au_m\|_H\right].
\]
Using (\ref{ustrong}) and a Gronwall inequality we obtain
\be
u_m\in L^{\infty}(0,T; {\mathcal{D}}(A)) \cap L^2(0,T; {\mathcal {D}}(A^{\fr{3}{2}}))
\la{super}
\ee
with uniform bounds in $m$, independent of $T$ and depending only on initial data. It is only now that we attained by interpolation the uniform bounds for $u_m\in L^2(0,T, B(\Omega))$.

Finally, we take the equation (\ref{pmqeq}), apply $\l$ and obtain, for $g =\l q$
\be
\pa_t g + u_m\cdot\na g + \l g + [\l, u_m\cdot \na ]q = 0
\la{geq}
\ee
We multiply by $4g^3$ and integrate, using Proposition \ref{cordoba} and H\"{o}lder inequalities
\[
\fr{d}{dt}\|g\|_{L^4(\Omega)}^4 \le 4\|g\|_{L^4(\Omega)}^3\|[u_m\cdot\na, \l q]\|_{L^4(\Omega)}
\]
Now, by (\ref{lfour})
\[
\|[u_m\cdot\na, \l q]\|_{L^4(\Omega)} \le C\|[u_m\cdot\na, \l]q\|_{\fr{1}{2},D}
\]
and, by (\ref{comtwo})
\[
\|[u_m\cdot\na, \l]q\|_{\fr{1}{2},D} \le C\|u_m\|_{B(\Omega)}\|q\|_{\fr{3}{2},D}\in L^1((0,T))
\]
This belongs uniformly to $L^1((0,T))$ because of (\ref{qhone}) and (\ref{super}).
We obtain thus
\be
q\in L^{\infty}(0,T; W^{1,4}(\Omega)
\la{qwonefour}
\ee
with bounds that are indpendent of $T$ and depend only on initial data. 

This concludes the uniform bounds, which are (\ref{qhone}), (\ref{super}) and (\ref{qwonefour}). Passage to the limit $m\to\infty$ is done using an Aubin-Lions lemma \cite{lions} and the bounds are inherited by the solutions of the limit equations. We omit further details.

\section*{Acknowledgements}
The work of PC was partially supported by NSF grant DMS-1209394. The work of VV was partially supported by NSF grant DMS-1514771 and by an Alfred P. Sloan Research Fellowship. The work of TE was partially supported by NSF grant DMS-1402357.

\end{document}